\documentclass[reqno,12pt]{amsart}
\newtheorem{theo}{Theorem}
\newtheorem{rem}{Remark}

\newtheorem{lem}{Lemma}

\newtheorem{df}{Definition}

\newcommand\eps\varepsilon
\newcommand\ph\varphi
\newcommand\kap\Lambda

\usepackage{enumerate}
\usepackage{graphicx}
\usepackage{float}

\begin{document}

\title[Bounded Solutions to Differential Inclusions]
{Forward-Backward Bounded Solutions to Differential Inclusions }

\author[Oleg Zubelevich]{Oleg Zubelevich\\ \\\tt
 \\Dept. of Theoretical mechanics,  \\
Mechanics and Mathematics Faculty,\\
M. V. Lomonosov Moscow State University\\
Russia, 119899, Moscow,  MGU \\ozubel@yandex.ru
 }
\email{ozubel@yandex.ru}
\date{}
\thanks{Partially supported by RFBR  18-01-00887.}

\begin{abstract}It is shown that some class of differential inclusions has solutions that are defined and bounded for all real values of independent variable. Applications to dynamics are considered.
\end{abstract}

\maketitle
\numberwithin{equation}{section}
\newtheorem{theorem}{Theorem}[section]
\newtheorem{lemma}[theorem]{Lemma}
\newtheorem{definition}{Definition}[section]

\section{Introduction}

 Discontinuous differential equations often arise in applied mathematics: one of the traditional examples is the motion
of a body in presence of dry friction \cite{filipp}; more recent sources of interest come from control theory and games
theory. In  important article \cite{aaa18} the author  studies and compares  solutions in sense of \cite{filipp}
with other notions, due to Krasowskii and Hermes, and with the classical ones (Newton and
Carathéodory solutions).  Other remarkable work was done in \cite{aaa6}, \cite{aaa11} and \cite{aaa8} (see also the reference therein).

In the context of control theory, other types of
solutions (Euler solutions) have been successfully employed (see Ancona et al \cite{bbb1}, Malisoff
et al \cite{bbb8}).

In the present article we focus our attention on discontinuous differential equations of type which  usually arises in mechanics of systems with the  Coulomb (dry) friction. Such systems are described by Filippov's construction presented in \cite{filipp}.

We study the systems of the Newton second law type, that is the systems of the second order equations.  The Coulomb  friction generates discontinuities in velocities and generally it does not spoil dependence on time and spatial variables.
So that we impose conditions which carry out such a specialization.  

We show that in some sense  unstable system with Coulomb friction has solutions defined for all real time and these 
solutions are bounded in both directions of time. Herewith the Coulomb friction mollifies the system in the direction of positive time and destabilizes it in the negative time direction.  

The main effect can be described as follows. Assume we have a mechanical system
$$\ddot x=-\frac{\partial V}{\partial x}(x),\quad x\in\mathbb{R}^m$$ 
and the potential energy $V=V(x)$ attains maximum at a point $\tilde x\in\mathbb{R}^m$.  It is well known that the equilibrium $x(t)\equiv \tilde x$ is unstable. Now let us perturbate  the system:
\begin{equation}\label{dtddddy55}\ddot x=-\frac{\partial V}{\partial x}(x)+g(t,x,\dot x).\end{equation}
It turns out that if we impose  certain conditions on $g,\quad V$ then   system  (\ref{dtddddy55}) has a solution $x(t)$ that is defined for all $t\in\mathbb{R}$ and $\sup_{t\in\mathbb{R}}|x(t)|<\infty$. This holds true even if the term $g$ is a discontinuous function of $\dot x$, the Coulomb friction for example.

\section{ Definitions and The Statement of the Problem}\label{secc1}

Let $\mu$ stand for the standard Lebesgue measure in $$\mathbb{R}^m=\{x=(x^1,\ldots, x^m)\},\quad d\mu_x=dx^1\ldots dx^m.$$
By $(x,y)=\sum_{k=1}^mx^iy^i$ we denote the standard Euclidean inner product and $|x|$ stands for the corresponding norm.

Let $$B_r(x_0)=\{x\in\mathbb{R}^m\mid |x-x_0|<r\}$$ stand for the open ball of radius $r>0$ with center at $x_0$;
$$\partial B_r(x_0)=\{x\in\mathbb{R}^m\mid |x-x_0|=r\}.$$

Let $\mathrm{conv}\, U,\quad U\subset\mathbb{R}^m$ stand for the closed convex hull of $U$. 

Let $M\subset \mathbb{R}^m$ be an open domain. Introduce a domain  $$G=\mathbb{R}\times M\times\mathbb{R}^m.$$
Consider a mapping $$ f:G\to \mathbb{R}^m,\quad f=f(t,x,y).$$

For each $(t,x)\in\mathbb{R}\times M$ the mapping $f(t,x,\cdot)$
is measurable in $y$;
for almost all $y$ the functions
$$f(t,x,y),\quad d_x f(t,x,y),\quad d^2_{xx}f(t,x,y)$$
are continuous  in $(t,x)\in\mathbb{R}\times M$.

Moreover, the following hypotheses hold:
\begin{description}
\item[A]
for any compact set $K\subset\mathbb{R}\times M$ there exists a positive constant $c_K$ such that for almost all
$y$  and for all $(t,x)\in K$ it follows that 
$$|d_x f(t,x,y)|+|d ^2_{xx}f(t,x,y)|+|f(t,x,y)|<c_K;$$
 \item[B] for all small enough $\eps>0$, for any compact interval $I\subset\mathbb{R}$ and for any compact set $K\subset M$ there exists $\delta>0$ such that for almost all $y$ and for all $$x\in K,\quad
t',t''\in I$$
the following implication holds:
$$|t'- t''|<\delta\Longrightarrow
|f(t',x,y)-f( t'', x,y)|<\eps;$$
\item[C] 
for all small enough $\eps>0$, for any $t$ and for any $\tilde x\in M$ there exists $\delta>0$ such that for almost all $y$ we have 
$$ f(t,B_\delta(\tilde x),y)\subset B_\eps(f(t,\tilde x,y)).$$
\end{description}

The main object of our study is the following initial value problem
\begin{equation}\label{xdfdd}
\ddot x=f(t,x,\dot x).\end{equation}
 Now we introduce a concept of generalized solution to this system. For briefness we will say "generalized solution" but it would be more accurate to call it "the solution in the sense of inclusions".

 The following definition  is a modified version of one in \cite{filipp}. This modification is physically reasonable.

\begin{df}\label{dfgff4} We shall say that a function $x(t)\in C^1((t_1,t_2),\mathbb{R}^m)$ is a generalized solution to problem
(\ref{xdfdd}) if

1)  $\dot x(t)$ is an absolutely continuous function in $(t_1,t_2)$;

2) for almost all $t\in(t_1,t_2) $ the following inclusion holds
\begin{equation}\label{dfg5}\ddot x(t)\in\bigcap_{r>0}\bigcap_{N}\mathrm{conv}\,f\Big(t,x(t),B_r\big(\dot x(t)\big)\backslash N\Big).\end{equation}
Here $\bigcap_{N}$ stands for the intersection over all measure-null sets: $$N\subset\mathbb{R}^m,\quad \mu(N)=0.$$\end{df}Recall that an absolutely continuous function has derivative almost everywhere and this derivative is locally Lebesgue integrable \cite{KF}.
\begin{rem}
If the function $f$ is continuous in $G$ then the  set in the right-hand side of (\ref{dfg5})  consists of the single element $\{f(t,x(t),\dot x(t))\}$.
\end{rem}

\section{The Main Theorem}\label{secc2}Introduce a function $F\in C^4(M)$ and two sets
$$D_c=\{x\in M\mid F(x)<c\},\quad \hat D_c=\{x\in M\mid F(x)=c\}.$$ Here $c$ is a constant.

Suppose that
\begin{description}
\item[d1]   the closure of $D_c$ in $\mathbb{R}^m$ is a compact set and it is contained in $M$:
$$\overline {D_c}\subset M;$$
\item[d2] there exists a homeomorphism $\psi:\overline {B_1(0)}\to D_c\cup \hat D_c$ such that
$$\psi(\partial B_1(0))=\hat D_c.$$
\end{description}
\begin{theo}\label{sdfddr}
Assume that {\bf A}, {\bf B}, {\bf C}, {\bf d1}, {\bf d2} are satisfied. Let a form $d^2F(x)$ 
be   positive definite or positive semi-definite for all $x\in \hat D_c$;
 and
for almost all $$(t,x,y)\in \mathbb{R}\times \hat  D_c\times\mathbb{R}^m$$ it follows that
$$
dF(x)[f(t,x,y)]>0.$$
Then equation (\ref{xdfdd}) has a generalized solution $x(t)\in C(\mathbb{R},\mathbb{R}^m)$, and 
for all $t$ one has
$$x(t)\in \overline D_c.$$  \end{theo}
\begin{rem}The theorem remains valid if we replace condition {\bf d2} with the following one:
the set $D_c$ does not admit a continuous retraction $\overline D_c\to \partial D_c$ and $\partial D_c=\hat D_c.$

If $f\in C^2(G,\mathbb{R}^m)$ and all the above conditions except {\bf C} are satisfied  then the theorem remains valid and  $x(t)$ is a solution in the classical sense.
\end{rem}
Note also that the function $\dot x(t)$ is not obliged  to be bounded in $\mathbb{R}$.

Theorem \ref{sdfddr} is proved in section \ref{qwe2}; section \ref{qwe1} contains auxiliary facts.

\section{Applications to Dynamics}\label{sdrgggg}
In this section we consider a pair of model examples. 
\subsection{System with Potential Forces and Discontinuous Perturbation}
Let $F,f$ be the same functions as in sections \ref{secc1}, \ref{secc2} and all the above conditions, particularly {\bf A}, {\bf B}, {\bf C}, {\bf d1}, {\bf d2},  are fulfilled.
 
 Consider a system
\begin{equation}\label{dfg556}\ddot x=\frac{\partial F}{\partial x}(x)+f(t,x,\dot x).\end{equation}
Here the function $F$ plays a role of potential energy taken with opposite sign. 
\begin{theo}\label{sdfffdr}
Let the form $d^2F(x)$ 
be   positive definite or positive semi-definite for all $x\in \hat D_c$;
 and
for almost all $$(t,x,y)\in \mathbb{R}\times \hat  D_c\times\mathbb{R}^m$$ it follows that
$$
|dF(x)|^2+dF(x)[f(t,x,y)]>0.$$
Then equation (\ref{dfg556}) has a generalized solution $x(t)\in C(\mathbb{R},\mathbb{R}^m)$, and 
for all $t$ one has
$$x(t)\in \overline D_c.$$  \end{theo}
Indeed, introduce a function $\eta\in C^\infty (\mathbb{R}^m)$ such that  
$\eta (x)=1$ provided $x\in \overline{D_c}$ and $\mathrm{supp}\,\eta\subset M$ is a compact set.

Consider a system 
$$\ddot x=\eta(x)\frac{\partial F}{\partial x}(x)+f(t,x,\dot x).$$
By theorem \ref{sdfddr} this system has a solution $x(t)$ and this solution belongs to $\overline{D_c}$ for all time. Thus $x(t)$ is a solution to system (\ref{dfg556}) also.

\subsection {Never Falling Pendulum}

{\it A pendulum  consists of a weightless rod $OA,\quad |OA|=R$ pivoted to a fixed ideal hinge  $O$. On the  end $A$ the rod has a particle of mass $m$. 
There is also a wheel of radius $r,\quad r<R$ pivoted by its center at the same point $O$. 

All the system is situated in a fixed vertical plane and experiences the standard gravity $\boldsymbol g$. 

The wheel rotates  about a horizontal  axis passing though the point $O$.  The  angular velocity of the wheel  $\omega=\omega(t)\in C^1(\mathbb{R})$ is a given function. There is a force of Coulomb friction between the wheel and the rod. The magnitude of this force is characterized by a coefficient $k>0$.

We assume that the wheel is not homogeneous so that the coefficient $k$ depends on mutual position of the rod and the wheel. 

This means that $K$ depends on $\psi$ (see below).  } 
\begin{figure}[H]
\centering
\includegraphics[width=45mm, height=42 mm]{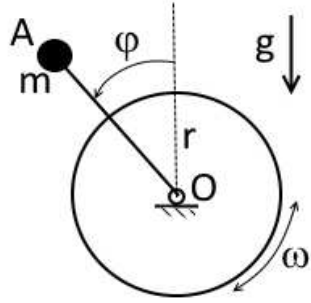}
\caption{The pendulum and the wheel. \label{oow}}
\end{figure}

Choose the dimensions such that
$$R=1,\quad m=1,\quad g=1.$$ The equation of motion is
$$\ddot\ph=\sin\ph-K\frac{\dot\ph-\omega}{|\dot\ph-\omega|},\quad K=kr.$$ The angle $\ph$ is shown at the picture.

After a change of variable
$$\ph=\psi+\Phi(t),\quad \Phi(t)=\int_{\mathrm{const}}^t\omega(\xi)d\xi$$ this equation takes the form
\begin{equation}\label{dty55}
\ddot \psi=\sin\big(\psi +\Phi(t)\big)-K\frac{\dot\psi}{|\dot\psi|}-\dot\omega(t).\end{equation}
Here $$K=K(\psi)\in C^2(\mathbb{R}),\quad K(\psi)\le K^*$$ is a non-negative  valued and $2\pi$-periodic function.

Definition \ref{dfgff4} implies that if   $\dot\psi=0$ for all $t$ from an open interval then
$$|\sin\big(\psi +\Phi(t)\big)-\dot\omega(t)|\le K(\psi).$$
This meets the Coulomb model of friction.

Assume that there are constants $\Phi^*,\eps$ such that
$$\sup_{t\in\mathbb{R}}|\Phi(t)|\le \Phi^*,\quad \sup_{t\in\mathbb{R}}|\dot\omega(t)|\le \eps.$$
\begin{theo}\label{dg555}
Suppose there exists a constant $\psi^*\ge 0$ such that  $$\psi^*+\Phi^*\le \pi/2,\quad\sin(\psi^*-\Phi^*)>K^*+\eps.$$

Then problem (\ref{dty55}) has a generalized solution $\psi\in C(\mathbb{R})$ such that 
$$|\psi(t)|\le\psi^*,\quad t\in\mathbb{R}.$$\end{theo}
Particularly this means that 
$$|\ph(t)|\le \psi^*+\Phi^*,\quad t\in\mathbb{R}$$ and 
 the pendulum never falls down. 

Theorem \ref{dg555} follows from theorem \ref{sdfddr} by taking
$$M=\mathbb{R},\quad F(\psi)=\psi^2,\quad D_c=(-\psi^*,\psi^*).$$

\section{ Regular Statement }\label{qwe1}
\subsection{Forward Bounded Solutions}
We use the following smooth result. 

Let
$$H=\mathbb{R}_+\times M\times\mathbb{R}^m,\quad \mathbb{R}_+=\{t\ge 0\}.$$
Consider a system
\begin{equation}\label{dfbbbbbg55}
\ddot x=a(t,x,\dot x),\quad a\in C^2(H,\mathbb{R}^m).\end{equation}

\begin{theo}[\cite{zubel}]\label{fggr5}
Assume that the conditions   {\bf d1}, {\bf d2} hold and 
for any compact set $$K\subset \mathbb{R}_+\times M$$ there exists a positive constant $C_K$ such that for  all $$(t,x,y)\in K\times \mathbb{R}^m$$ one has
$$|a(t,x,y)|<C_K.$$

Assume also that for all $x\in \hat D_c$ a quadratic form $d^2F(x)$ is  positive definite or positive semi-definite; and
for all $$(t,x,y)\in \mathbb{R}_+\times \hat  D_c\times\mathbb{R}^m$$ it follows that
$$
dF(x)[a(t,x,y)]>0.$$
Then equation (\ref{dfbbbbbg55}) has a  solution $x(t)\in C^4(\mathbb{R}_+,\mathbb{R}^m)$, and 
for all $t\ge 0$ one has
$$x(t)\in \overline D_c.$$  \end{theo}

\subsection{Forward-Backward bounded solutions}
Consider a system
\begin{equation}\label{bbbbg55}
\ddot x=b(t,x,\dot x).\end{equation}
The functions $b=b(t,x,y)$,
$$d_{y}b(t,x,y),\quad d_{x}b(t,x,y),\quad d^2_{xy}b(t,x,y),\quad d^2_{yy}b(t,x,y),\quad d^2_{xx}b(t,x,y)$$ are continuous in $G$.

Let the following hypotheses hold. 

 \begin{description}
\item[H1]
For any compact set $$K\subset\mathbb{R}\times M$$ there exists positive constant $\check c_K$ such that for  all $$(t,x,y)\in K\times \mathbb{R}^m$$ one has
\begin{align}
|d_{x}b(t,x,y)|&+| d^2_{xy}b(t,x,y)|+| d^2_{yy}b(t,x,y)|+|d^2_{xx}b(t,x,y)|\nonumber\\&+
|b(t,x,y)|+|d_{y}b(t,x,y)|<\check c_K;\nonumber\end{align}
 \item[H2] 
for all small enough $\eps>0$, for any compact interval $I\subset\mathbb{R}$ and for any compact set $K\subset M$ there exists $\delta>0$ such that for  all $(x,y)\in  K\times\mathbb{R}^m$ and for all
$t',t''\in I$
the following implication holds:
$$|t'- t''|<\delta\Longrightarrow
|b(t',x,y)-b( t'', x,y)|<\eps.$$
\end{description}

 \begin{theo}\label{fgqqwegr5}
Assume that the conditions {\bf H1}, {\bf H2}, {\bf d1}, {\bf d2} hold.

Assume also that for all $x\in \hat D_c$ a quadratic form $d^2F(x)$ is  positive definite or positive semi-definite; and
for all $$(t,x,y)\in \mathbb{R}\times \hat  D_c\times\mathbb{R}^m$$ it follows that
$$
dF(x)[b(t,x,y)]>0.$$
Then equation (\ref{bbbbg55}) has a  solution $x(t)\in C^2(\mathbb{R},\mathbb{R}^m)$ such that
$$x(t)\in \overline D_c$$ for any $t\in\mathbb{R}.$ \end{theo}

\subsubsection{ Proof of theorem \ref{fgqqwegr5}} Let $\chi_k(t),\quad k\in \mathbb{N}$ stand for the indicator:
$$\chi_k(t)=1\quad\mbox{provided}\quad t\in(-k,k],$$
and $\chi_k(t)=0$ otherwise.

Pick a function $\psi\in C^\infty(\mathbb{R})$ such that $$\mathrm{supp}\,\psi\subset (-1,1),\quad \psi\ge 0,\quad\int_{\mathbb{R}}\psi(t)dt=1.$$
Introduce a sequence of functions
$$\delta_l(t)=l\psi( lt),\quad \mathrm{supp}\,\delta_l\subset (-1/l,1/l),\quad l\in\mathbb{N}.$$

Introduce a function 
$$b_k(t,x,y)=\sum_{j\in\mathbb{Z}}b(t+2kj,x,y)\chi_k(t+2kj)$$
and a function
$$ b_{kl}(t,x,y)=\int_{\mathbb{R}}b_k(s,x,y)\delta_l(s-t)ds.$$
It is clear all the functions $b_k,b_{kl}$ are $2k-$ periodic in $t$ and 
\begin{equation}\label{dfhgr6u6u}b_k(t,x,y)=b(t,x,y),\quad t\in (-k,k].\end{equation}
Moreover $b_{kl}\in C^2(G,\mathbb{R}^m).$

For any $\eps>0$, for any $T>0$ and for any compact set $K\subset M$ there exists $L>0$ such that for all $k\ge 2T$,
$$t\in[-T,T]\subset[-k/2,k/2],\quad(x,y)\in K\times\mathbb{R}^m$$ the following implication holds
\begin{equation}\label{dhd5y4y}l>L\Longrightarrow|b_{kl}(t,x,y)-b(t,x,y)|<\eps.\end{equation}
This follows from the hypothesis {\bf H2} and formula (\ref{dfhgr6u6u}).

Consider the following system
\begin{equation}\label{dfgrrrr}
\ddot x=b_{kl}(t,x,\dot x).\end{equation}
This system satisfies all the conditions of theorem \ref{fggr5}.
Indeed,
let $(t,x,y)\in \mathbb{R}\times \hat  D_c\times\mathbb{R}^m$ then
\begin{align}
dF&(x)[b_{kl}(t,x,y)]\nonumber\\
&=\int_\mathbb{R}\sum_{j\in\mathbb{Z}}dF(x)[b(s+2kj,x,y)]\chi_k(s+2kj)\delta_l(s-t)ds>0;\nonumber\end{align}
 and for any compact set  $K\subset\mathbb{R}\times M$ it follows that
 \begin{equation}\label{sdfg5689}|b_{kl}(t,x,y)|\le\check  c_{ K},\quad (t,x,y)\in K\times \mathbb{R}^m.\end{equation}
 Thus system (\ref{dfgrrrr}) has a solution $x_{kl}\in C^4(\mathbb{R}_+,\mathbb{R}^m)$ and 
 $$x_{kl}(t)\in \overline {D_c},\quad t\ge 0.$$
 
 Introduce  functions 
$$z_{kli}(t)=x_{kl}(t+2ki),\quad i\in\mathbb{N}.$$
 All these functions are the solutions to (\ref{dfgrrrr}) and 
 \begin{equation}\label{dfggggg46}z_{kli}\in C^4([-2ki,\infty),\mathbb{R}^m);\qquad z_{kli}(t)\in\overline {D_c},\quad t\ge -2ki.\end{equation}

 Due to (\ref{sdfg5689}) for $k,i,l\in\mathbb{N}$ and fixed $T\le k/2$ we have
 \begin{equation}\label{sdffff}
|\ddot z_{kli}(t)|\le \check c_{[-T,T]\times\overline{D_c}},\quad t\in[-T,T] .\end{equation}
 Let $H^{2}([-T,T],\mathbb{R}^m)$ stand for the Sobolev space of functions $u(t)$ such that 
$$u,\dot u,\ddot u\in L^2([-T,T],\mathbb{R}^m).$$
Recall that $H^2([-T,T],\mathbb{R}^m)$ is a Hilbert space with the inner product
$$(a,b)_{H^2[-T,T]}=\int_{-T}^T\big((a(s),b(s))+(\ddot a(s),\ddot b(s))\big)ds.$$
Another inner product 
$$\langle a,b\rangle_{H^2[-T,T]}=\int_{-T}^T\big((a(s),b(s))+(\dot a(s),\dot b(s))+(\ddot a(s),\ddot b(s))\big)ds$$
gives the same topology in $H^2([-T,T],\mathbb{R}^m)$.

Recall also that there is a compact embedding \begin{equation}\label{dfdf456}H^2([-T,T],\mathbb{R}^m)\subset C^1([-T,T],\mathbb{R}^m).\end{equation}
Here and below we refer \cite{adams} for the properties of the Sobolev spaces.

\begin{lem}\label{sfg44}The sequence $Z_p:=z_{ppp}$ contains a subsequence $Z_{p_s}$ that is convergent 
to a function $z_*\in C^1(\mathbb{R},\mathbb{R}^m)$ in $C^1([-T,T],\mathbb{R}^m)$ for any $T>0$.

The function $z_*$ is such that $$z_*(t)\in \overline{D_c},\quad t\in\mathbb{R}.$$\end{lem}

 {\it Proof of the lemma \ref{sfg44}.} 
 Fix $T=1$. From formulas (\ref{sdffff}),  (\ref{dfggggg46}), (\ref{dfdf456}) it follows that there exists a subsequence $Z_{p_q}$ that is convergent in $C^1([-1,1],\mathbb{R}^m)$. By the same reason this subsequence contains a subsequence that is convergent in $C^1([-2,2],\mathbb{R}^m)$ etc.

The diagonal argument finishes the proof.

The function $z_*$ is the announced solution  to system (\ref{bbbbg55}). 

Indeed, to see this fix $T>0$ and for $t\in[-T,T]$ write 
$$\dot Z_{p_s}(t)=\dot Z_{p_s}(0)+\int_0^tb_{p_sp_s}(\xi,Z_{p_s}(\xi),\dot Z_{p_s}(\xi))d\xi.$$
To pass to the limit in this equality  let
\begin{align}
\int_0^t&b_{p_sp_s}(\xi,Z_{p_s}(\xi),\dot Z_{p_s}(\xi))d\xi-\int_0^tb(\xi,z_*(\xi),\dot z_*(\xi))d\xi\nonumber\\
&=\int_0^t\big(b_{p_sp_s}(\xi,Z_{p_s}(\xi),\dot Z_{p_s}(\xi))-b(\xi,Z_{p_s}(\xi),\dot Z_{p_s}(\xi))\big)d\xi\nonumber\\
&+\int_0^t\big(b(\xi,Z_{p_s}(\xi),\dot Z_{p_s}(\xi))-b(\xi,z_*(\xi),\dot z_*(\xi))\big)d\xi.\label{xdfg990}
\end{align}
By (\ref{dhd5y4y}) we have
$$|b_{p_sp_s}(\xi,Z_{p_s}(\xi),\dot Z_{p_s}(\xi))-b(\xi,Z_{p_s}(\xi),\dot Z_{p_s}(\xi))|\to 0$$ uniformly in $\xi\in[-T,T]$.

On the other hand the difference 
$$|b(\xi,Z_{p_s}(\xi),\dot Z_{p_s}(\xi))-b(\xi,z_*(\xi),\dot z_*(\xi))|\to 0$$ vanishes
pointwise in $\xi\in[-T,T]$ and 
$$|b(\xi,Z_{p_s}(\xi),\dot Z_{p_s}(\xi))-b(\xi,z_*(\xi),\dot z_*(\xi))|\le2 \check c_{[-T,T]\times\overline{D_c}}.$$ Thus the last integral in (\ref{xdfg990})  also tends to zero by the dominated convergence theorem.

Theorem \ref{fgqqwegr5} is proved.

\section{Proof of Theorem \ref{sdfddr}}\label{qwe2}

Pick a function $\ph\in C^\infty(\mathbb{R}^m)$ such that $$\mathrm{supp}\,\ph\subset B_1(0),\quad \ph\ge 0,\quad\int_{\mathbb{R}^m}\ph(x)d\mu_x=1.$$
Introduce a sequence of functions
$$\delta_k(x)=k^m\ph( kx),\quad \mathrm{supp}\,\delta_k\subset B_{1/k}(0),\quad k\in\mathbb{N}$$
and put
\begin{equation}\label{sdfs55}f_k(t,x,y)=\int_{\mathbb{R}^m}f(t,x,z)\delta_k(z-y)d\mu_z.\end{equation} 

Constructions of type (\ref{sdfs55}) are usually employed in approximation theory. Intuitively speaking, the sequence $f_k$ approximates  the function $f$ in some sense. Nevertheless, in this section we do not use such an argument and  we do not refer  any approximation theorems.

We use only one convolution property: $f_k(t,\cdot,\cdot)\in C^2(M\times\mathbb{R}^m,\mathbb{R}^m).$

Observe that the following system
\begin{equation}\label{sfddet}
\ddot x=f_k(t,x,\dot x)\end{equation} satisfies all the conditions of theorem \ref{fgqqwegr5}.
Indeed, {\bf H2} follows from {\bf B};
\begin{equation}\label{df778}|f_k(t,x,y)|\le\int_{\mathbb{R}^m}c_K\delta_k(z-y)d\mu_z\le c_K;\end{equation}
and
$$dF(x)[f_k(t,x,y)]=\int_{\mathbb{R}^m}dF(x)[f(t,x,z)]\delta_k(z-y)d\mu_z>0$$ provided 
$(t,x)\in\mathbb{R}\times\hat D_c.$

Thus theorem \ref{fgqqwegr5} supplies each  system (\ref{sfddet}) with the solution $x_k(t)$. 

From the properties of the Lebesgue integral for  all $t$ we have 
\begin{align}
\ddot x_k(t)&=f_k(t,x_k(t),\dot x_k(t))\nonumber\\
&\in\bigcap_N\mathrm{conv}\,f\Big(t,x_k(t),B_{\frac{1}{k}}\big(\dot x_k(t)\big)\backslash N\Big).\label{fsgff}\end{align}
Here the intersection is taken over all measure-null sets:
$$N\subset\mathbb{R}^m,\quad \mu(N)=0.$$
Indeed,
\begin{align}f_k&(t,x_k(t),\dot x_k(t))=\int_{\mathbb{R}^m}f(t,x_k(t),z)\delta_k(z-\dot x_k(t))d\mu_z\nonumber\\
&=
\int_{\mathbb{R}^m\backslash N}f(t,x_k(t),z)\delta_k(z-\dot x_k(t))d\mu_z\nonumber\\
&\in
\mathrm{conv}\,f\Big(t,x_k(t),B_{\frac{1}{k}}\big(\dot x_k(t)\big)\backslash N\Big).\nonumber\end{align}

\begin{lem}\label{xsssdfg55}
For any $T>0$ the sequence $\{x_k\}$ is bounded in $$H^2([-T,T],\mathbb{R}^m),$$ that is
$$\sup_k\|x_k\|_{H^2[-T,T]}<\infty.$$\end{lem}
{\it Proof of lemma \ref{xsssdfg55}.} First, notice that $$\{x_k(t)\}\subset \overline{D_c}$$ for all $t\in\mathbb{R}$ and $k$. 
Thus by  formula (\ref{df778})
we get
$$|\ddot x_k(t)|\le c_{[-T,T]\times\overline {D_c}},\quad t\in[-T,T] .$$  The lemma is proved.

\begin{lem}\label{dfg5687}
The sequence $\{x_k\}$ contains a subsequence $\{x_{k_s}\}$ that is convergent to a function 
$$x_*\in H^2_{\mathrm{loc}}(\mathbb{R},\mathbb{R}^m)$$ in the following sense:

for any $T>0$ one has
$$\|x_{k_s}-x_*\|_{C^1[-T,T]}\to 0,$$
and $x_{k_s}\to x_*$ weakly in $H^{2}([-T,T],\mathbb{R}^m)$.
\end{lem}

{\it Proof of lemma \ref{dfg5687}.}
Recall a theorem. \begin{theo}[\cite{iosida}]\label{cdfgh6789} A Banach space  $(Y,\|\cdot\|_Y)$ is reflexive iff any bounded sequence in $Y$ contains a weakly convergent subsequence.
\end{theo}

Take an increasing sequence $T_n\to\infty.$

Since the embedding $H^2([-T_n,T_n],\mathbb{R}^m)\subset C^1([-T_n,T_n],\mathbb{R}^m)$ is compact and the space $H^2([-T_n,T_n],\mathbb{R}^m)$ is reflexive,
we can choose a subsequence $\{x_{k_j}^{(1)}\}\subset \{x_k\}$ that is convergent in $C^1([-T_1,T_1],\mathbb{R}^m)$ and weakly convergent in $$H^2([-T_1,T_1],\mathbb{R}^m).$$ The sequence $\{x_{k_j}^{(1)}\}$ contains a subsequence $\{x_{k_j}^{(2)}\}$ that is convergent in 
$C^1([-T_2,T_2],\mathbb{R}^m)$ and weakly convergent in $H^2([-T_2,T_2],\mathbb{R}^m)$ etc. The diagonal sequence is convergent in the desired manner.

The lemma is proved.

Note that $x_*(t)\in \overline{D_c},\quad t\in\mathbb{R}$.

Choose a constant $\rho>0$ such that for any $x_0\in \overline{D_c}$ it follows that $\overline{B_r(x_0)}\subset M,\quad r\le\rho.$

Introduce sets
$$U(t,r,r')=\bigcap_{\mu(N)=0}\mathrm{conv}\,f\Big(t,B_r\big(x_*(t)\big),B_{r'}\big(\dot x_*(t)\big)\backslash N\Big)\subset\mathbb{R}^m,\quad r<\rho.$$
By lemma \ref{ddrrr} (see below) the sets $U(t,r,r')$ are nonvoid for almost all $t$.
The sets $U(t,r,r')$ are  closed and convex as an intersection of closed convex  sets.

The sets $U(t,r,r')$ are uniformly bounded relative $$t\in [-T,T]\quad ,0<r\le\rho,\quad r'>0.$$ Indded,
\begin{equation}\label{dhgt}z\in U(t,r,r')\Longrightarrow |z|\le c_{[-T,T]\times S_r},\quad S_r=\overline{\bigcup_{x\in D_c}B_r(x)}\subset M. \end{equation}

Let $W^T(r,r')\subset L^2([-T,T],\mathbb{R}^m)$ stand for a set of functions $u(t)\in L^2([-T,T],\mathbb{R}^m)$ such that for almost all $t\in[-T,T]$ one has
$$u(t)\in U(t,r,r').$$

\begin{lem}\label{dfgfffg45}The sets $W^T(r,r')$ are closed bounded and convex.\end{lem}
{\it Proof of lemma \ref{dfgfffg45}.} Convexity is evident. Prove that the sets are closed. Indeed, let a sequence $\{w_i\}\subset W^T(r,r')$ converges to $w$ in $L^2([-T,T],\mathbb{R}^m)$. Then it contains a subsequence $\{w_{i_j}\}$ that is convergent to $w$ almost everywhere \cite{folland}. Thus $w(t)\in U(t,r,r')$ for almost all $t$.

The boundedness follows from  formula (\ref{dhgt}).

The lemma is proved.

\begin{lem}\label{ddrrr}For any $T,r'>0,\quad r\in(0,\rho]$ there exists a number $J$ such that
$$k_s>J\Longrightarrow \ddot x_{k_s}\in W^T(r,r').$$\end{lem}
{\it Proof of lemma \ref{ddrrr}.} Choose $J$ such that

1) $J>2/r'$;

2) $k_s>J\Longrightarrow \|x_{k_s}-x_*\|_{C^1[-T,T]}<\min\{r,r'/2\}.$

Take any $k_s>J$ it is clear
\begin{equation}\label{dfgdff5}x_{k_s}(t)\in B_r(x_*(t)).\end{equation}
The following inclusion holds for all $t\in[-T,T]$
\begin{equation}\label{xdfggg445t}B_{\frac{1}{k_s}}(\dot x_{k_s}(t))\subset B_{r'}(\dot x_*(t)).\end{equation}
Indeed,
Let $\xi\in B_{\frac{1}{k_s}}(\dot x_{k_s}(t))$ then
$$|\xi-\dot x_*(t)|\le |\xi-\dot x_{k_s}(t)|+|\dot x_{k_s}(t)-\dot x_*(t)|<\frac{1}{k_s}+\frac{r'}{2}<r'.$$
The assertion of the lemma follows from formulas (\ref{dfgdff5}), (\ref{xdfggg445t}), (\ref{fsgff}).

The lemma is proved.

\begin{lem}\label{ddrdfgdfggrr}For any $T,r'>0,\quad r\in(0,\rho]$ 
one has
$$\ddot x_*\in W^T(r,r').$$\end{lem}
Indeed, by lemma \ref{dfg5687} the sequence $\{\ddot x_{k_s}\}$ is weakly convergent to $\ddot x_*$ in $L^2([-T,T],\mathbb{R}^m)$. The set  $W^T(r,r')$ is convex and closed thus it is weakly closed \cite{edvards}. This proves the lemma.

The result of lemma \ref{ddrdfgdfggrr} can explicitly be formulated as follows. For almost all $t$
we have
\begin{equation}\label{sdg567}\ddot x_*(t)\in \bigcap_{r>0}\bigcap_{r'>0}\,\bigcap_{\mu(N)=0}\mathrm{conv}\,f\Big(t,B_r\big(x_*(t)\big),B_{r'}\big(\dot x_*(t)\big)\backslash N\Big).\end{equation}

According to {\bf C} for any $t$ we can choose  sequences 
$$\delta_q,\eps_q\to 0,\quad \mathbb{N}\ni q\to\infty$$ and a set $\mathcal N,\quad \mu(\mathcal N)=0$
such that for all $q$ and for all $y\in\mathbb{R}^m\backslash\mathcal N$ the following inclusion holds
\begin{equation}\label{dasdrg55}f\big(t,B_{\delta_q}(x_*(t)),y\big)\subset B_{\eps_q}\big(f(t,x_*(t),y)\big).\end{equation}

\begin{lem}\label{dfg5768}For any
$r'>0$ and for any \begin{equation}\label{dfg667}N,\quad \mathcal N\subset N,\quad \mu(N)=0\end{equation}
the following equality holds
\begin{align}\bigcap_{q\in\mathbb{N}}&\mathrm{conv}\,f\Big(t,B_{\delta_q}\big(x_*(t)\big),B_{r'}\big(\dot x_*(t)\big)\backslash  N\Big)
\nonumber\\&=\mathrm{conv}\,f\Big(t,x_*(t),B_{r'}\big(\dot x_*(t)\big)\backslash  N\Big).\label{dgh66h}\end{align}\end{lem}
Obviously, formula (\ref{dgh66h}) can be rewritten as follows
\begin{align}\bigcap_{r>0}&\mathrm{conv}\,f\Big(t,B_{r}\big(x_*(t)\big),B_{r'}\big(\dot x_*(t)\big)\backslash  N\Big)
\nonumber\\&=\mathrm{conv}\,f\Big(t,x_*(t),B_{r'}\big(\dot x_*(t)\big)\backslash  N\Big).\label{dfg56jj}\end{align}
Moreover,
observe that 
\begin{align}\bigcap_{\mu(N)=0}\,\bigcap_{r>0}&\mathrm{conv}\,f\Big(t,B_{r}\big(x_*(t)\big),B_{r'}\big(\dot x_*(t)\big)\backslash N\Big)
\nonumber\\&=
\widetilde\bigcap\,\bigcap_{r>0}\mathrm{conv}\,f\Big(t,B_{r}\big(x_*(t)\big),B_{r'}\big(\dot x_*(t)\big)\backslash N\Big),\nonumber
\end{align}
and 
\begin{align}\bigcap_{\mu(N)=0}&\mathrm{conv}\,f\Big(t,x_*(t),B_{r'}\big(\dot x_*(t)\big)\backslash N\Big)
\nonumber\\&=
\widetilde\bigcap\mathrm{conv}\,f\Big(t,x_*(t),B_{r'}\big(\dot x_*(t)\big)\backslash N\Big),\nonumber
\end{align}

where $\bigcap_{\mu(N)=0}$ means the intersection over all measure-null sets $N$ and $\widetilde\bigcap$ means the intersection over all $N$ such that (\ref{dfg667}).

So that formula (\ref{dfg56jj}) implies
\begin{align}\bigcap_{\mu(N)=0}\bigcap_{r>0}&\mathrm{conv}\,f\Big(t,B_{r}\big(x_*(t)\big),B_{r'}\big(\dot x_*(t)\big)\backslash  N\Big)
\nonumber\\&=\bigcap_{\mu(N)=0}\mathrm{conv}\,f\Big(t,x_*(t),B_{r'}\big(\dot x_*(t)\big)\backslash  N\Big).\nonumber\end{align}
Thus formula (\ref{sdg567}) takes the form 
$$\ddot x_*(t)\in \bigcap_{r'>0}\,\bigcap_{\mu(N)=0}\mathrm{conv}\,f\Big(t,x_*(t),B_{r'}\big(\dot x_*(t)\big)\backslash N\Big).$$
This proves  theorem \ref{sdfddr}.

{\it Proof of lemma \ref{dfg5768}.} 
From formula (\ref{dasdrg55})
for   all $y\in B_{r'}\big(\dot x_*(t)\big)\backslash  N$ we have 
$$ f\Big(t,B_{\delta_q}\big(x_*(t)\big),y\Big)\subset B_{\eps_q}\Big(f\big(t,x_*(t),y\big)\Big).$$

So that the set 
\begin{align}
f&\Big(t,B_{\delta_q}\big(x_*(t)\big),B_{r'}\big(\dot x_*(t)\big)\backslash N\Big)\nonumber\\
&=\bigcup_{y\in B_{r'}\big(\dot x_*(t)\big)\backslash N}f\Big(t,B_{\delta_q}\big(x_*(t)\big),y\Big)\nonumber\\
&\subset\bigcup_{y\in B_{r'}\big(\dot x_*(t)\big)\backslash N} B_{\eps_q}\Big(f\big(t,x_*(t),y\big)\Big)\nonumber\end{align}

is contained in an $\eps_q-$ neighbourhood of the  set 
$$f\Big(t,x_*(t),B_{r'}\big(\dot x_*(t)\big)\backslash N\Big).$$
Thus we obtain
\begin{equation}\label{dfh68yu}\bigcap_{q\in\mathbb{N}}f\Big(t,B_{\delta_q}\big(x_*(t)\big),B_{r'}\big(\dot x_*(t)\big)\backslash  N\Big)=
\overline{f\Big(t,x_*(t),B_{r'}\big(\dot x_*(t)\big)\backslash N\Big)}.\end{equation}
Observe that the closed convex hull of a set $Q\subset\mathbb{R}^m$ is the intersection over all closed half-spaces that contain $Q$. That is
$$\mathrm{conv}\,Q=\bigcap_{\xi\in\mathbb{R}^m}P_\xi,\quad P_\xi=\{u\in\mathbb{R}^m\mid(u,\xi)\le \sup_{x\in Q}(x,\xi)\}.$$
Observe also that \begin{equation}\label{dfgg6}\mathrm{conv}\,\overline{Q}=\mathrm{conv}\,Q.\end{equation}

We have
\begin{equation}\label{dfg66}
\bigcap_{q\in\mathbb{N}}\mathrm{conv}\,f\Big(t,B_{\delta_q}\big(x_*(t)\big),B_{r'}\big(\dot x_*(t)\big)\backslash N\Big)=
\bigcap_{q\in\mathbb{N}}\bigcap_{{\xi\in\mathbb{R}}^m}P_{q\xi},\end{equation} where
\begin{align}P_{q\xi}&=\Big\{u\in\mathbb{R}^m\Big| (u,\xi)\nonumber\\&\le\sup\big\{(x,\xi)\mid x\in f\big(t,B_{\delta_q}\big(x_*(t)\big),B_{r'}\big(\dot x_*(t)\big)\backslash  N\big)\big\}\Big\}\nonumber.\end{align}
Changing the order of intersections in the right-hand side of formula (\ref{dfg66}) and by formula (\ref{dfh68yu}) we get 

\begin{align}\bigcap_{q\in\mathbb{N}}&\mathrm{conv}\,f\Big(t,B_{\delta_q}\big(x_*(t)\big),B_{r'}\big(\dot x_*(t)\big)\backslash  N\Big)
\nonumber\\&=\mathrm{conv}\,\overline{f\Big(t,x_*(t),B_{r'}\big(\dot x_*(t)\big)\backslash N\Big)}.\nonumber\end{align}
Now formula (\ref{dgh66h}) follows from (\ref{dfgg6}).

The lemma is proved.


\end{document}